\documentclass[12pt,leqno,fleqn]{article}
\usepackage{amssymb, amsmath, amsthm}
\usepackage{mathrsfs}       
\usepackage{color}  
     
\newcommand{\R}{\mathbb{R}}

\newcommand{\N}{\mathbb{N}}

\newcommand{\esssup}{\operatorname*{ess\,sup}}

\newcommand{\trace}{\operatorname{trace}}

\newcommand{\supp}{\operatorname*{supp}}

\newcommand{\dist}{\operatorname*{dist}}

\newcommand{\const}{\operatorname*{const}}

\newcommand{\dsum}{\displaystyle\sum}

\newcommand{\bq}{\begin{eqnarray}}
\newcommand{\eq}{\end{eqnarray}}
\newcommand{\bqn}{\begin{eqnarray*}}
\newcommand{\eqn}{\end{eqnarray*}}

\newcommand{\intl}{\int\limits}

\newcommand{\Beweisende}{\rule{0.2cm}{0.2cm}}

\newcommand{\D}{\displaystyle}

\newcounter{secnum}

\newtheorem{thm}{Theorem}[section]

\newtheorem{cor}[thm]{Corollary}
\newtheorem{lem}[thm]{Lemma}

\theoremstyle{definition}

\newtheorem{defin}[thm]{Definition}
\newtheorem{rem}[thm]{Remark}

\title{On the Serrin-type condition on  one  velocity component for  the  Navier-Stokes  equations
} %Liouville problem for the steady 3D Navier-Stokes equations}
 
\author{Dongho Chae$^*$  and J\"{o}rg Wolf $^\dagger$\\
\ \\
Department of Mathematics\\
Chung-Ang University\\
 Seoul 06974, Republic of Korea\\
 ($*$)e-mail: dchae@cau.ac.kr\\
($\dagger$)e-mail: jwolf2603@cau.ac.kr}
\date{}
\begin{document}
\maketitle
\begin{abstract}
In this paper we consider the regularity problem of the Navier-Stokes equations in $ \R^{3} $. We show that the 
Serrin-type condition imposed on one  component  of the velocity $ u_3\in L^p(0,T; L^q(\R^{3} ))$    
 satisfying  $ \frac{2}{p}+ \frac{3}{q} <1$, $ 3<q \le +\infty$  
implies the regularity of  the weak Leray  solution $ u: \R^{3} \times (0,T) \rightarrow \R^{3} $ with the initial data
belonging to $ L^2(\Bbb R^3)   \cap L^3(\R^{3})$. The result is an immediate consequence of a new local regularity criterion in terms of one velocity component 
for suitable weak solutions. 
\noindent{\bf keywords:} Navier-Stokes equation, regularity of weak solutions, Serrin condition for one velocity component

\end{abstract}
%\tableofcontents

%%% ----------------------------------------------------------------------
%       SECTION 1
%%% ----------------------------------------------------------------------
\section{Introduction}
\label{sec:-1}
\setcounter{secnum}{\value{section} \setcounter{equation}{0}
\renewcommand{\theequation}{\mbox{\arabic{secnum}.\arabic{equation}}}}

 Let $ 0< T< +\infty$ and let $ \Omega  \subset \R^{3} $ be a domain. We consider the Navier-Stokes equations in the space time cylinder $ \Omega _T= \Omega  \times  (0,T)$
\begin{equation}
\begin{cases}
\partial _t u+ (u\cdot \nabla) u -\Delta  u= -\nabla \pi,
\\[0.2cm]
\nabla \cdot u=0,
\end{cases}
\label{nse}
\end{equation}
equipped with the initial and boundary condition 
\begin{align}
u &=  u_0  \quad  \text{ on}\quad  \Omega \times \{ 0\},  
 \label{initial}
\\
u &=0\quad   \text{on}\quad  \partial \Omega \times  (0,T), 
\label{bound}
\end{align}
where $u=(u_1, u_2, u_3)=u(x,t)$ represents  the velocity of the fluid flows, and $\pi =\pi (x,t)$ denotes the scalar pressure. 
The existence of global  weak solutions has been proved by Leray \cite{ler} and Hopf \cite{hopf}. However, the existence 
of global regular  solutions remains as  an outstanding open problem.  In order to get  deeper  understanding for the regularity problem 
of the Navier-Stokes equations people  study  various sufficient conditions which guarantee the regularity. 
The first of such  conditions was introduced independently by  Prodi \cite{prod} and Serrin \cite{serrin}, namely, 
if the weak solution $ u$ satisfies the condition 
\begin{equation}
 u\in L^p(0,T; L^q(\Omega )),\quad  \frac{2}{p} + \frac{3}{q} \le 1,\quad 2 \le p < +\infty, 
\label{serrin}
 \end{equation}  
 then   $ u$ is regular. The case $ p=+\infty$ and $ q=3$ has been proven  later  by  Escauriaza,  Seregin  and \v{S}ver\'ak in \cite{essesv}. A similar condition 
also holds for local domains.  Clearly, condition  \eqref{serrin}  with equality case is invariant under  the following natural scaling of the Navier-Stokes equations 
\begin{equation}
u(x,t)\mapsto u^\lambda (x,t)=\lambda u(\lambda x, \lambda^2 t),\quad  \lambda >0.   
\label{scal}
 \end{equation}

Concerning the partial regularity, it  has been proved by Caffarelli,Kohn and Nirenberg in 
\cite{caf} that suitable weak solutions, which satisfy the local energy inequality (cf.  \cite{sch1}, \cite{sch2}, \cite{caf}), the one-dimensional parabolic Hausdorff measure of the set of possible singularities is zero.  This   result is a consequence of  the 
$ \varepsilon-$regularity criterion imposed  on quantities invariant under the scaling   \eqref{scal}. 
A new class of sufficient conditions have been introduced in  \cite{nepe}, \cite{nenope}  imposed   only on one component, say $  u_3\in L^p(0,T; L^q(\R^{3} ))$ with $ \frac{2}{p}+ \frac{3}{q} \le \frac{1}{2}$, 
which is much stronger than  \eqref{serrin}. This condition was improved by several authors up to a condition of the form
\begin{equation}
 u_3\in L^p(0,T; L^q( \R^{3} )),\quad  \frac{2}{p} + \frac{3}{q} \le  1-\delta (q),\quad \frac{10}{3} \le  q < +\infty,
\label{serrin1}
 \end{equation}  
where $ \delta (q)= \frac{1}{2}- \frac{1}{2q}$ (cf. \cite{cz1}, \cite{pz1}).  
Other conditions imposed on reduced components of the gradient of velocity,  which are invariant under the natural scaling  \eqref{scal} 
have been established by various authors  (cf. \cite{kz1}, \cite{cz1}, \cite{cz2}, \cite{cati1}, \cite{cati2}, \cite{wol7}). 
For conditions on the pressure see \cite{ss1}, \cite{zh}, \cite{cl1}, \cite{beirsga}, \cite{beir2} and on the vorticity 
see \cite{beir1}, \cite{bers}, \cite{cc} (for more discussion on  this  topic see \cite{nest}). 
Furthermore,  in \cite{bc} it has been proved that the Serrin condition \eqref{serrin} imposed on two components guarantees the regularity (for the localization of this result see \cite{baewolf1}). 

\vspace{0.2cm}
The aim of the present paper is to remove the number $ \delta (q)$ 
in  \eqref{serrin1} and get the regularity under the almost Serrin condition (Serrin's condition omitting equality)
\begin{equation}
 u_3\in L^p(0,T; L^q( \R^{3}  )),\quad  \frac{2}{p} + \frac{3}{q} <1,\quad 3 <   q \le  +\infty. 
\label{serrin1}
 \end{equation}

\vspace{0.2cm}
Throughout the paper we use the following function spaces. By $ W^{m,\,q}(\Omega ), W^{m,\,q}_0(\Omega ), m\in \N, 1 \le q \le +\infty$ we denote the usual Sobolev spaces. 
If no confusion arises,  we use the same notation  for  spaces of vector valued functions.   
The space $ L^p_{ \sigma }(\Omega )$ stands for the closure of  the space of smooth solenoidal functions with compact support   
$ C^\infty_{ c, \sigma }(\Omega )$   with respect to  the $ L^p$ norm.  In addition, for functions $ f\in L^p_\sigma(\Omega )$ 
it holds   $ f\cdot n=0$ on $\partial  \Omega $ in the sense of distributions, whenever 
$ f\in L^p_{ \sigma }(\Omega )$ (for more details see \cite{sohr}).  

Let $ -\infty \le a < b \le +\infty$. Given a normed space $ X$ by $ L^p(a,b; X), 1 \le p \le +\infty$, we denote the 
space of Bochner measurable functions $ f: (a, b) \rightarrow X$, such that 
\[
\begin{cases}
{\D \|f\|_{ L^p(a,b; X)} = \bigg(\intl_{a}^b \|f(s)\|^p ds\bigg) ^{ \frac{1 }{p}}< +\infty}\quad   &\text{if}\quad  1 \le p < +\infty,
\\[0.3cm]
\|f\|_{ L^\infty(a,b; X)}=  {\D \esssup_{s\in (a,b)}} \|f(s)\| &\text{if}\quad   p= +\infty. 
\end{cases}
\] 
By $ V^{1, 2}(\Omega _T)$, we denote the energy space 
$ L^\infty(0,T; L^2(\Omega  ))\cap L^2(0,T; W^{1,\,2}(\Omega  ) )$. In addition, we define the following sub spaces of $ V^{1,2}(\Omega _T)$,
\begin{align*}
 V^{ 1,2}_0(\Omega _T) &= V^{ 1,2} (\Omega _T)\cap  L^2(0,T; W^{1,\,2}_0(\Omega  ) ),
\\
V^{ 1,2}_\sigma (\Omega _T) &=  V^{ 1,2} (\Omega _T)\cap L^\infty(0,T; L^2_\sigma (\Omega  )),
\\
V^{ 1,2}_{ 0, \sigma }(\Omega _T) &=  V^{ 1,2}_\sigma(\Omega _T) \cap  L^2(0,T; W^{1,\,2}_0(\Omega  ) ). 
\end{align*}
Note that by virtue  of Sobolev's inequality  and H\"older's inequality we have for all $ u\in V^{ 1,2}(\Omega _T)$, 
\begin{equation}
 u \in L^s(0,T; L^q(\R^{3} )),\quad  \forall 2 \le q  \le +\infty, \quad  \frac{2}{s}+ \frac{3}{q} = \frac{3}{2}. 
\label{1.4}
 \end{equation} 

Given matrices $ A, B \in \R^{3\times  3} $ by $ A:B$ we denote the scalar product $  \sum_{i,j=1}^{3} A_{ ij}B_{ ij} = 
\trace(AB^{ \top})$. Then $ |A|=  (A:A)^{ \frac{1}{2}} $ denotes the Euclidian  norm in $ \R^{3\times  3} $. For 
two vectors $ a, b\in \R^{3} $ by $ a\cdot b$ and $ |a|$ we denote the usual scalar product and norm in $ \R^{3} $ respectively. 
For $ x_0\in \R^{3} $ and $ 0<r< +\infty$ we denote  $ B(x_0, r)$ the ball in $ \R^{3} $  with center $ x_0$ and 
radius $ r$. Given $ x_0 \in   \R^3  $    and $ 0<\rho, r<+\infty$, we define  the non isotropic spatial  cylinder 
\[
 U(x_0;\rho, r):= 
   B'(x'_0, \rho )\times  (x_{ 0,3}- r, x_{ 0,3}+r),
\]
where 
  \[
  B'(x'_0,\rho )= \{(x_1, x_2)\in \R^{2}\,|\, (x_1-x_{ 0,1})^2 + (x_2-x_{ 0,2})^2 < \rho ^ 2\},\quad  x'_0=(x_{ 0,1}, x_{ 0,2}).
\]
 For $ z_0= (x_0, t_0)\in \R^{3} \times  \R$ 
we denote by $ Q(z_0; \rho, r)$ the  non isotropic parabolic cylinder  $ U(x_0; \rho, r)\times  (t_0-r^2, t_0)$.

\vspace{0.2cm}
We now recall the notion of a suitable weak solution to  \eqref{nse}. 

\begin{defin}
 \label{def1.1}
 A pair $ (u, \pi )\in V^{1, 2}_{\sigma }(\Omega_T)\times  L^{ \frac{3}{2}}(\Omega _T)$ is called a {\it suitable weak solution} 
 to  \eqref{nse} if $ u$ solves  \eqref{nse} in the sense of distributions, and satisfies the local energy inequality
  \begin{align}
&\frac{1}{2} \intl_{ \Omega  } |u(t)|^2 \phi (t) dx + \intl_{0}^t\intl_{\Omega  } |\nabla u|^2 \phi  dx ds   
\cr
 &\le  \frac{1}{2}\intl_{0}^t\intl_{\Omega} |u|^2 (\partial _t +  \Delta ) \phi  dx ds   +
\frac{1}{2} \intl_{0}^t\intl_{ \Omega } (|u|^2 + 2\pi  )  u\cdot  \nabla \phi dx ds
\label{en}
\end{align} 
for all non negative $ \phi \in C^\infty_c(\Omega \times  (0,T] )$ and for almost  all $ t\in (0,T)$.
\end{defin}

Our  main result concerns the  regularity conditions of almost Serrin-type imposed on  one velocity component in case of the whole space $ \Omega = \R^{3}$. 
\begin{thm}
 \label{thm1.4}
Let $ u_0    \in L^2_\sigma (\R^{3} )\cap L^3(\R^{3} )$. 
Let $ u\in V^{ 1,2}_{ 0,\sigma }(\R^{3} \times  (0, T))$ be a weak Leray solution 
 to  \eqref{nse},  \eqref{initial}. 
   Suppose  the following almost Serrin condition for one velocity component holds 
 \begin{equation}
 \begin{cases}
 u_3 \in L^p(0, T; L^q( \R^{3}))
 \\[0.3cm]
  \text{for some}\quad  3 <  q \le  +\infty\quad   \text{with}\quad  \frac{2}{p}+ \frac{3}{q} <1.
 \end{cases}
 \label{serrin-w}
 \end{equation}  
Then, $ u$ is a regular solution. 

\end{thm} 

In fact, the proof of Theorem\,\ref{thm1.4} is based on the following more general  result of local regularity as a consequence 
of a new local regularity criterion imposed on one velocity component.    

\begin{thm}
 \label{thm1.1}
 Let $ (u, \pi )\in V^{ 1,2}_{\sigma }(\Omega _T)\times  L^{ \frac{3}{2}}(\Omega _T)$ be a suitable weak solution 
 to  \eqref{nse}. Let $ z_0= (x_0, t_0) = (x_{ 0}', x_{03}, t_0 )\in \Omega _T$.  Assume, $ u_3\in L^p(t_0-\rho ^2, t_0; L^q( U(x_0; \rho, \rho )))$ 
 for some $ 0< \rho  < \min\{\dist(x_0, \partial \Omega ),  \sqrt{T}\}$ and $ 3 \le q < +\infty$, with $ 
 \frac{2}{p} + \frac{3}{q} \in \Big[ 1, \frac{3}{2}\Big]$.  

  1. Suppose, 
 \begin{equation}
  \limsup_{ r \searrow 0} r^{1-  \frac{2}{p} - \frac{3}{q}} \| u_3\|_{ L^p(t_0 - r^2, t_0; 
  L^q(U(x_0;\rho , r ))} =0.
 \label{cond1}
  \end{equation} 
  Then for all $ 0< \lambda  <1$ it holds
  \begin{align}
&\sup_{ 0<r<\rho }  r^{ -\lambda }\Big( \| u\|^2_{ L^\infty(t_0 - r^{ 2}, t_0; L^2(U(x_0; \rho , r))} 
   +  \| \nabla u\|^2_{ L^2(t_0 - r^{ 2}, t_0; L^2(U(x_0;\rho , r))}\Big)< +\infty.  
  \label{morrey}
  \end{align}
  2. Suppose furthermore that for some $ 1 < \alpha   < +\infty$, 
  \begin{equation}
  \limsup_{ r \searrow 0}  (-\log (r))^{ \alpha } r ^{1-  \frac{2}{p} - \frac{3}{q}} \| u_3\|_{ L^p(t_0 - r^2, t_0; L^q
  (U(x_0; \rho , r))} =0.
 \label{cond2}
  \end{equation} 
 Then $ z_0$ is a regular point.  
  \end{thm}

\begin{rem}
 \label{rem1.2a}
Clearly, \eqref{cond1}  is satisfied if $ u_3$ satisfies the 
 Serrin condition  with $ \frac{2}{p} + \frac{3}{q} \le1$, and  the first part of the theorem says  in this case that 
  \begin{equation}
 u \in  L^\infty_{ loc}((0,T]; {\cal M}^{ 2, \lambda }_{ loc}(\Omega ))$ and $ \nabla u \in {\cal M}^{ 2, \lambda }_{ loc}(\Omega \times  (0,T])\quad \forall 0< \lambda <1,
  \label{cond2a}
  \end{equation} 
 where $ {\cal M}_{ loc}^{2, \lambda }(\Omega )$ stand for the local Morrey space. However,  \eqref{cond2a}  is not 
 sufficient to  guarantee the  regularity of $u$.  On the other hand, as we will see in Section\,3, if those Morrey conditions hold with $ \lambda =1$ the 
 regularity follows. This is the reason why we need to add the logarithmic factor in  \eqref{cond2} in order to get the regularity.    
 Also we wish to remark that condition  \eqref{cond2} defines  the  set regular points in terms of $ u_3$ 
 which leads to the  partial regularity in terms of one velocity component.

\end{rem}

As a corollary of Theorem\,\ref{thm1.1} we get the following almost Serrin local condition in terms of one velocity component.
\begin{cor}
 \label{cor1.3}
Let $ (u, \pi )\in V^{ 1,2}_{\sigma }(\Omega _T)\times  L^{ \frac{3}{2}}(\Omega _T)$ be a suitable weak solution 
 to  \eqref{nse}. Let $ Q(z_0, R)  \subset \Omega _T$, where $ z_0=(x_0, t_0)$. Suppose that the following  condition holds 
 \begin{equation}
 \begin{cases}
 u_3 \in L^p(t_0, t_0-R^2; L^q(B(x_0, R)))
 \\[0.3cm]
  \text{for}\quad  3 <  q \le  +\infty \quad   \text{with}\quad  \frac{2}{p}+ \frac{3}{q} <1.
 \end{cases}
 \label{serrin-y}
 \end{equation}  
Then, all points in $ B(x_0, r)\times  (t_0-\rho ^2, t_0]$ are regular points.
\end{cor}

\begin{rem}
 \label{rem1.3} The assumption on existence of the pressure in Theorem\,\ref{thm1.1} and Corollary\,\ref{cor1.3} 
 is not essential, and we can replace the notion of suitable weak solution by local suitable weak solution, using the local pressure projection introduced  in \cite{wol8}.   
\end{rem}

%%% ----------------------------------------------------------------------
%       SECTION 
%%% ----------------------------------------------------------------------
\section{Proof of Theorem\,\ref{thm1.4} and Corollary\,\ref{cor1.3}}
\label{sec:-5}
\setcounter{secnum}{\value{section} \setcounter{equation}{0}
\renewcommand{\theequation}{\mbox{\arabic{secnum}.\arabic{equation}}}}

{\bf Proof of Theorem\,\ref{thm1.4}}:
Since $ u_0  \in L^3(\R^{3})$  by local well posedness of the Navier-Stokes equations in $L^3(\R^{3} )$, 
there exists  maximal time $ 0< T_{ \ast} \le +\infty$ such that $ u\in C([0;T_{ \ast}); L^3(\R^{3} ))$ (see e.g. Kato \cite{kato}). Assume $ 0< T_{ \ast} \le T$ 
and it holds $ \|u(t)\|_{ L^3(\R^{3} )} \rightarrow +\infty$ as $ t \nearrow T_{ \ast}$. Let $ z_0= (x_0, t_0)$ with $ t_0= T_{ \ast}$ be fixed. Since $ u$ is regular in $ \R^{3} \times  (0, T_{ \ast})$ there exists $ \pi \in L^{ 3/2}(\R^{3} \times  (0,T))$ 
and $ (u, \pi )$ is a suitable weak solution. There exists $ 2 \le \widetilde{p} < p \le +\infty$ such that 
\[
\frac{2}{\widetilde{p}} + \frac{3}{q} =1.  
\]
Thus by means of H\"older's inequality we get for $ 0< r \le \rho $ with  $ \rho =  \sqrt{T_{ \ast}}$
\begin{align*}
\|u_3\|_{ L^{ \widetilde{p}}(t_0- r^2, t_0; L^{ q}(U(x_0, \rho , r)))}
&\le
cr^{2 \frac{p-\widetilde{p}}{p\widetilde{p}}}\|u_3\|_{ L^{p}(t_0- r^2, t_0; L^{ q}(U(x_0, \rho , r)))} 
\\
&\le cr^{2 \frac{p-\widetilde{p}}{p\widetilde{p}}} \|u_3\|_{ L^p(0,T; L^q(\R^{3} ))}   
\end{align*}
with a constant $ c>0$ independent of $ r$. Thus, since $(-\log(r))^{ \alpha } r^{2 \frac{p-\widetilde{p}}{p\widetilde{p}}}
\rightarrow 0$ as $ r \rightarrow 0$ for any $ \alpha >1$ 
condition  \eqref{cond2} of Theorem\,\ref{thm1.1} is fulfilled. This implies that $ (x_0, t_0)= (x_0, T_{ \ast})$ 
is a regular point. In particular, $ u\in L^\infty(0, T_{ \ast}; L^3(\R^{3} ))$ which contradicts to the definition of $ T_{ \ast}$. 
Accordingly, the assertion of the Theorem is true.  
 \hfill \Beweisende
 
 \vspace{0.3cm}
 {\bf Proof of Corollary\,\ref{cor1.3}}: Set $ z_0=(x_0, t_0)$. Let $ (y_0, s_0) \in B(x_0, R)\times  (t_0-R ^2, t_0]$. There exists $ 0< \rho <R$ 
 such that $ U(y_0, \rho, \rho  )\times  (s_0-\rho ^2, s_0)  \subset Q(z_0, R)$. Arguing as in the proof of Theroem\,\ref{thm1.4},  we find    $ 2 \le \widetilde{p} < p \le +\infty$ such that 
$  \frac{2}{\widetilde{p}} + \frac{3}{q} =1$, and thus by means of H\"older's inequality 
 \begin{align*}
\|u_3\|_{ L^{ \widetilde{p}}(s_0- r^2, s_0; L^{ q}(U(y_0,\rho, r)))}
\le r^{2 \frac{p-\widetilde{p}}{p\widetilde{p}}} \|u_3\|_{ L^p(t_0- R^2,t_0; L^q(B(x_0, R)))}.    
\end{align*} 
Thus, condition  \eqref{cond2} of Theorem\,\ref{thm1.1} is fulfilled, which  implies that $ (y_0, s_0)$ 
is a regular point.  \hfill \Beweisende

%%% ----------------------------------------------------------------------
%       SECTION 
%%% ----------------------------------------------------------------------
\section{Proof of Theorem\,\ref{thm1.1}}
\label{sec:-2}
\setcounter{secnum}{\value{section} \setcounter{equation}{0}
\renewcommand{\theequation}{\mbox{\arabic{secnum}.\arabic{equation}}}}

 In the proof of Theorem\,\ref{thm1.1} we make essential use of the following 
\begin{lem}
 \label{lem3.1}
Let $ \frac{1}{2} \le R \le 1$ and $ 0< r< 1$. Set $ U= B'(R)\times  (-r, r)$ and  $ Q=  U\times  (-r^2,0)$. 
Then  $ L^\infty(-r^2, 0; L^2(U))\cap L ^2(-r^2, 0; W^{1,\,2} (U))$ is continuously embedded into   
$L^m(-r^2,0; L^l(U))$ for all 
$ 2 \le m \le \infty, 2 \le l \le 6$ such that $ \frac{m}{2}+ \frac{3}{l} = \frac{3}{2}$ and the following inequality holds 
\begin{align}
\|u\|_{ L^m(-r^2,0; L^l(U))}^2 \le c \Big(\|u\|^2_{ L^\infty(-r^2,0; L^2(U))} + \|\nabla u\|^2_{ L^2(Q)}\Big), 
\label{emb1}
\end{align}  
where $ c>0$ stands for  an absolute  constant.  
\end{lem}

{\bf Proof}:  1. First we show that 
\begin{equation}
 \|u\|_{ L^6(U)} \le cr^{ -1} \| u \|_{ L^2(U)} + c\|\nabla u\|_{ L^2(U)},\quad  \forall   u\in W^{1,\,2} (U), 
 \label{emb2}
 \end{equation} 
where  $ c>0$ stands for an  absolute constant.  In fact,  given  $ u\in W^{1,\,2}(U) $, setting  $ \widetilde{U}= B'(R)\times  (-2r, 2r)$, we define  the extension 
\[
\widetilde{u}(x):=\begin{cases}
u(x)\quad   &\text{if}\quad  x\in U,
\\[0.3cm]
u(x', 2r- x_3) &\text{if}\quad  x\in \widetilde{U}, x_3\in [r, 2r), 
\\[0.3cm]
u(x', -2r- x_3) &\text{if}\quad  x\in \widetilde{U}, x_3\in (-2r, -r].
\end{cases}
\] 

Then $ \widetilde{u}\in W^{1,\,2} (\widetilde{U})$ and it holds
\begin{equation}
 \|\nabla \widetilde{u}\|_{ L^2(\widetilde{U})} \le 3 \|\nabla u\|_{ L^2(U)},\quad  
 \|\widetilde{u}\|_{ L^2(\widetilde{U})} \le 3 \|u\|_{ L^2(U)}. 
\label{ext1}
 \end{equation} 
Let $ \eta \in  C^\infty_c(-2r,2r)$ denote a cut off function  such that $ \eta =1$ on $ (-r,r)$ and $ |\eta '| \le 2r^{ -1} $. 
Noting that by $ r \le 2R$ it holds  $ \widetilde{U}  \subset B'(R) \times  (-4R, 4R) $ 
it holds   $ \widetilde{u} \eta \in W^{1,\,2} (B'(R)\times  (-4R, 4R))$.  By means of Sobolev's inequality  and a simple 
scaling argument along with  \eqref{ext1}, we get 
\begin{align*}
\| u\|_{ L ^6(U)} &\le \| \widetilde{u} \eta \|_{ L ^6(B'(R)\times  (-4R, 4R))}  
\\
&\le c\|\widetilde{u}\eta \|_{ L^2(B'(R)\times  (-4R, 4R))} + c\|\nabla (\widetilde{u}\eta )\|_{ L^2(B'(R)\times  (-4R, 4R))} 
\\
&\le cr^{ -1}\|\widetilde{u} \|_{ L^2(\widetilde{U}))} + c\|\nabla \widetilde{u}\|_{ L^2(\widetilde{U})} 
\\
&\le cr^{ -1}\|u \|_{ L^2(U))} +c\|\nabla u\|_{ L^2(U)}.
\end{align*}
Whence,  \eqref{emb1}

\vspace{0.3cm}
2. Let $ u\in L^\infty(-r^2, 0; L^2(U))\cap L ^2(-r^2, 0; W^{1,\,2} (U))$. 
Thanks to  \eqref{emb1} we have for almost all $ t\in (-r^2,0)$, 
\begin{align*}
\|u(t)\|^2_{ L^6(U)} \le c r^{ -2} \|u(t)\|^2_{ L^2(U)}   + c \|\nabla u(t)\|^2_{ L^2(U)},
\end{align*}
where $ c>0$ denotes an absolute constant. Integrating this inequality over $ (-r^2, 0)$ with respect to time, we arrive at
\begin{align*}
\|u\|^2_{ L^2(-r^2,0; L^6(U))} &\le cr^{ -2} \|u\|_{ L^2(-r^2,0; L^2(U))}^2 + c \|\nabla u\|^2_{ L^2(Q)}
\\
&\le c  \|u\|_{ L^\infty(-r^2,0; L^2(U))}^2  + c\|\nabla u\|^2_{ L^2(Q)}.
\end{align*}
Let $ 2 \le m \le \infty, 2 \le l \le 6$ such that $ \frac{m}{2}+ \frac{3}{l} = \frac{3}{2}$. 
By the aid of H\"older's inequality along with Young's inequality we find 
\begin{align*}
\|u\|^2_{ L^m(-r^2, 0; L^l(U))} &\le  \|u\|^{ \frac{4}{m}}_{ L^\infty(-r^2,0; L^2(U))}\|u\|^{ ^{ \frac{2m-4}{m}}}_{ L^2(-r^2,0; L^6(U))}
\\
&\le \|u\|^2_{ L^\infty(-r^2,0; L^2(U))} + \|u\|^2_{ L^2(-r^2,0; L^6(U))}.
\end{align*} 
Combining the last two inequalities, we obtain the desired estimate  \eqref{emb2}.  \hfill \Beweisende

\vspace{0.5cm}
{\bf Proof of Theorem\,\ref{thm1.1}.} By translation in space-time we may consider the case $z_0=(x_0, t_0)=(0,0)$ only.
First let us consider the case  that $ (u, \pi )\in V^{ 1,2}_\sigma (Q(0; 1,1))\times  L^{ \frac{3}{2}}(Q(0; 1,1))$ is a suitable weak solution to 
 \eqref{nse}.  Given $ 3< q< +\infty, 2< p < +\infty$ with $ \frac{2}{p} + \frac{3}{q} \in  \Big[1,\frac{3}{2}\Big]$, we denote 
 \[
Y(R) = \sup_{ 0< r \le R} r^{ 1- \frac{2}{p} - \frac{3}{q}} \| u_3\|_{ L^p(- r^2, 0; L^q(U(0; 1, r)))}.
\]

 Our aim will be to argue as in the proof of Caffarelli-Kohn-Nirenberg theorem  (cf. \cite{caf}), inserting  $ \phi =\Phi _{ n} 
 \zeta  $  in  \eqref{def1.1}, where $ \zeta $ denotes a cut-off function,  while  $ \Phi _n$ stands for the shifted fundamental solution to the backward heat equation in one spatial dimension, i.e.
\[
\Phi _n(x, t) = \frac{1}{  \sqrt{4\pi (-t+r_n^2)}} e^{ -\frac{x_3^2}{4(-t +r_n^2) }},\quad  (x,t)\in\R^3\times  (-\infty, 0),
\]   
where 
\[
r_n = 2^{ -n}, \quad  n\in \N. 
\]
In what follows we use the following notations. Given $ 0< R < +\infty$, we set 
\begin{align*}
U_n(R) &=U(0, R, r_n)=B'(R) \times  (-r_n, r_n), 
 \\
 Q_n(R) &= U_n \times  (-r_n^2, 0),
 \\
A_n(R) &= B'(R)\times  A^{ \ast}_n,
\end{align*}
where 
\begin{align*}
A^{ \ast}_n &= Q^{ \ast}_{ n} \setminus Q^{ \ast}_{ n+1},
\\
  Q^{ \ast}_n &=   (-r_n, r_n) \times  (-r_n^2, 0).
\end{align*}

In case $ R=1$ we write $Q_n  , A_n $ etc. in place of $ Q_n(1), A_n(1)$ etc.   

 Clearly, there exist absolute constants $ c_1, c_2>0$ such that  for all $ 0< R< +\infty$, $ n\in \N$ and $ j=1, \ldots, n$ it holds 
 \begin{align}
  & c_1 r_{ j}^{ -1} \le \Phi _n \le c_2 r_j^{ -1},\quad   c_1 r_{ j}^{ -2} \le |\partial _3\Phi _n| \le c_2 r_j^{ -2}\quad  
  \text{in}\quad A_j(R), 
 \label{2.4j}
 \\
 & c_1 r_{ n}^{ -1} \le \Phi _n \le c_2 r_n^{ -1},\quad   c_1 r_{ n}^{ -2} \le |\partial _3\Phi _n| \le c_2 r_n^{ -2}\quad  
  \text{in}\quad  Q_n(R).   
 \label{2.4n}
 \end{align}

Given  $ 0< R \le 1$ and $ n\in \N_0$,  the following  notation will be used in what follows 
 \begin{align*}
E_n(R) &= E_n^{ (1)}(R)+ E_n^{ (2)}(R):=  \esssup_{t\in (-r_n^2, 0)}\intl_{U_n(R)} |u(t)|^2   dx + 
\intl_{- r_n^2}^0 \intl_{U_n(R)} |\nabla u|^2   dx dt  , 
 \\
 {\cal E} &= E_0(1) =\esssup_{t\in (-1, 0)}\intl_{U_0(1)} |u(t)|^2  dx + 
\intl_{- 1}^0 \intl_{U_0(1)} |\nabla u|^2  dx dt  .
 \end{align*}

 In view of Lemma\,\ref{lem3.1},   we see that for all $ n\in \N_0$ and 
$ \frac{1}{2} \le  R \le 1$, 
 \begin{equation}
\| u  \|_{ L^m(-r_n^2, 0; L^l(U_n(R)))}^2 \le cE_n(R) \quad   \forall 2 \le m  \le +\infty, \quad  \frac{2}{m}+ \frac{3}{l} = \frac{3}{2}. 
 \label{2.6}
 \end{equation}

Let $ \eta= \eta (x_3, t) \in C^\infty_c((-1, 1)\times  (-1,0])$ denote a cut off function 
such that $ 0 \le \eta  \le 1$ in $ \R\times  (-1,0]$, $ \eta \equiv 1$ on $ Q_{ 1}^{ \ast}= (-1/2, 1/2)\times  (-1/4,0)$. 
In addition, let $ 1/2 \le \rho < R \le 1 $ be arbitrarily chosen, but fixed. 
Let $ \psi = \psi  (x') \in C^\infty ( \R^{2} )$,  such that $ 0 \le \psi  \le 1$ in $B'(R)$, $ \psi \equiv 0 $ 
in $ \R^{2} \setminus B( \frac{R+\rho }{2})$, 
 $ \eta \equiv 1$ on $ B'( \rho  )$, such that $ |\nabla ' \psi | \le c (R-\rho )^{ -1}, |D^2 \psi | \le  c (R-\rho )^{ -2}$, where $ \nabla ' = 
 (\partial _1, \partial _2)$ and $ x' = (x_1, x_2)$.
We insert $ \phi= \Phi _n \eta\psi$   into  \eqref{en}. This yields for $ n\in \N$,
 \begin{align}
&\frac{1}{2} \intl_{ U_0(R)} |u(\cdot, t)|^2 \Phi_n (\cdot , t) \eta(\cdot , t) \psi dx + \intl_{-1}^t\intl_{ U_0(R)} |\nabla u|^2 \Phi_n \eta \psi   dx ds   
\cr
 &\le  \frac{1}{2}\intl_{-1}^t\intl_{ U_0(R)} |u|^2 (\partial _t +  \Delta ) (\Phi_n  \eta\psi ) dx ds   
 +
\frac{1}{2} \intl_{-1}^t\intl_{ U_0(R) } |u|^2 u\cdot  \nabla  (\Phi_n \eta \psi )dx ds
\cr
&\qquad + \intl_{-1}^t\intl_{ U_0(R)} \pi   u \cdot \nabla  (\Phi_n \eta \psi ) dx ds =I+II+III.     
\label{en-bn}
\end{align} 
We now estimate $ I, II, III$ for $ n\in \N_0$. 
Recalling that $ (\partial _t + \partial _3\partial _3 )\Phi _n=0$ in $ Q_0(R)$,  and observing  \eqref{2.4j} and  \eqref{2.4n}, we get 
\begin{align*}
I &= \frac{1}{2}\intl_{-1}^t\intl_{ U_0(R)} |u|^2( \Phi _n \partial _t \eta\psi  + 2 \partial _3 \Phi _n \partial _3 
\eta\psi  + \Phi_n \Delta  (\eta \psi )  ) dx ds
\\
 &\le c  \intl_{A_0(R)} |u|^2 dxds + c(R-\rho )^{ -2}\intl_{Q_0(R)} |u|^2\Phi _n dxds 
 \\
 & \le  c(R-\rho )^{ -2}\sum_{j=0}^{n} r_j^{ -1}
 \intl_{Q_j(R)} |u|^2 dxds 
 \\
 &\le c(R-\rho )^{ -2}\sum_{j=0}^{n} r_j {\cal E} \le  c(R-\rho )^{ -2} {\cal E}. 
\end{align*}

Again using  \eqref{2.4j} and  \eqref{2.4n}, we find  
\begin{align*}
II & \le  \sum_{i=0}^{n-1} \intl_{A_i(R)} |u|^2 |u_3||\partial _3 \Phi _n| \eta \psi  dxds   + 
\intl_{Q_n(R)} |u|^2 |u_3| |\partial _3\Phi _n| \eta  \psi  dxds
\\
&\qquad \qquad  +
\intl_{Q_0(R)} |u|^3  \Phi _n \eta |\nabla' \psi  |  dxds
\\
& \le  c\sum_{j=0}^{n} r_j^{ -2}\intl_{Q_j(R)} |u|^2 |u_3|    dxds   + 
\intl_{Q_0(R)} |u|^3  \Phi _n  |\nabla' \psi  |dxds
\\
 & = II_1+ II_2.
\end{align*}
In our discussion below for given $ m\in (1, +\infty)$, by $ m'$ we denote the H\"older conjugate  $ \frac{m}{m-1}$. In case 
$ m=1$ we define  $ m'= +\infty$, while in case $ m=+\infty$ define  $ m'=1$.

We first estimate $ II_1$. Using H\"older's inequality, we find for $ j=0, \ldots, n$
 \begin{align*}
 r_j^{ -2}\intl_{Q_j(R)} |u|^2 |u_3|   dxds \le c Y(r_j) r_j^{ -3+ \frac{2}{p}+ \frac{3}{q}} \Bigg(\intl_{-r_j^2}^{ 0} \Bigg(    \intl_{U_j(R)} |u |^{2q'}  dx\Bigg)^{ \frac{p'}{q'}} ds\Bigg)^{ \frac{1}{p'}}. 
\end{align*}
Setting $ m= \frac{4q}{3}$ and $ l= 2q'$, we get $ \frac{2}{m} + \frac{3}{l}= \frac{3}{2}$. Applying H\"older's inequality again, 
we see that 
 \begin{align*}
 r_j^{ -2}\intl_{Q_j(R)} |u|^2 |u_3|  dxds    \le cY(r_j) r_j^{ -3 + \frac{2}{p}- \frac{3}{q} + \frac{2q-3p'}{qp'}} \Bigg(\intl_{-r_j^2}^{ 0} \Bigg(    \intl_{U_j(R)} |u |^{l}  dx\Bigg)^{ \frac{m}{l}} ds\Bigg)^{ \frac{2}{m}}.  
\end{align*}
Verifying that  $-3 + \frac{2}{p}- \frac{3}{q} + \frac{2q-3p'}{qp'}=-1$ and using   \eqref{2.6}, we obtain 
\[
  II_1 \le c  \sum_{j=0}^{n-1}r_j^{ -1} E_j(R)Y(r_j).
\]
Observing  \eqref{2.4j} and  \eqref{2.4n}, we estimate 
\[
II_2 \le c  (R-\rho )^{ -1}\sum_{j=0}^{n} r_j^{ -1}\intl_{Q_j(R)} |u|^3  dxds.
\]
Applying H\"older's inequality along with  \eqref{2.6}, we infer 
\begin{align*}
r_j^{ -1}\intl_{Q_j(R)} |u|^3 dxds \le 
r_j^{ - \frac{1}{2}} \|u\|^2_{ L^4(-r_j^2, 0; L^3(U_j(R))} {\cal E}^{ \frac{1}{2}}  \le c r_j^{ - \frac{1}{2}} 
E_j(R)  {\cal E}^{ \frac{1}{2}}.
\end{align*}
This gives
 \[
II_2 \le c(R-\rho )^{ -1}  {\cal E}^{ \frac{1}{2}}\sum_{j=0}^{n}r_j^{ - \frac{1}{2}} 
E_j(R).
\]

Accordingly,
\[
II \le  c  \sum_{j=0}^{n}  r_j^{ -1}E_j(R)Y(r_j)
+ c(R-\rho )^{ -1}  {\cal E}^{ \frac{1}{2}}\sum_{j=0}^{n}r_j^{ - \frac{1}{2}} 
E_j(R).
\]
It remains to estimate the integral  involving the  pressure. We write
\begin{align*}
III &= \intl_{-1}^t\intl_{U_0(R)} \pi u_3  \partial _3 \Phi _n \eta \psi   dxds + 
\intl_{-1}^t\intl_{U_0(R)}  \pi  \Phi _nu \cdot  \nabla (\eta \psi ) dxds= III_1+ III_2.       
\end{align*}
We define 
\[
\pi _0 = \mathscr{J}(u \otimes u \chi _{ U_0(R)} ), 
\]
where $ \mathscr{J}: L^m(Q_0(R)) \rightarrow L^m(\R^{3}\times  (-1,0) ), 1< m< +\infty$, stands for the bounded 
operator defined in the appendix, with the property
\[
-\Delta \pi _0 = \nabla \cdot \nabla \cdot (u \otimes u \chi_{ Q_0(R)}) \quad  \text{in}\quad \R^{3} \times  (-1,0), 
\]
in the sense of distributions. 
Thus, setting $ \pi _h= \pi - \pi _0$, it follows that $ \pi _h$ is harmonic in $ Q_0(R)$, and the following estimate holds true
\begin{align}
   \|\pi _h\|_{ L^{ \frac{3}{2}}(Q_0(R))} &\le c \|\pi\|_{ L^{ \frac{3}{2}}(Q_0(R))}
 +  \|\pi _0\|_{ L^{ \frac{3}{2}}( \R^{3} \times  (-1,0))}  
\cr
 &\le   c \|\pi\|_{ L^{ \frac{3}{2}}(Q_0(R))}
 +  c\|u\|_{ L^{3}(Q_0(R))}^2 \le c    (\|\pi\|_{ L^{ \frac{3}{2}}(Q_0(R))}+ {\cal E}).  
\label{hp}
\end{align}
Clearly, 
\[
III_1= \intl_{-1}^t\intl_{U_0(R)} \pi_0 u_3  \partial _3 \Phi _n \eta^2 \psi ^2 dxds + \intl_{-1}^t\intl_{U_0(R)} \pi_h u_3  \partial _3 \Phi _n \eta^2 \psi ^2 dxds = III_{ 11}+III_{ 12}.
\] 
Let $\sigma  \in (q ,+\infty)$ and $l \in [2, +\infty)$ such that 
\begin{equation}
 \frac{2}{\sigma } + \frac{3}{l} = 3. 
\label{2.50a}
 \end{equation} 
This implies that $ \frac{1}{l} = 1- \frac{2}{3\sigma  } > 
1-\frac{2}{3q} > \frac{1}{q'}$, i.e. $ l < q'$ and recalling $ \frac{2}{p} + \frac{3}{q} \le \frac{3}{2}$, we get 
$  \frac{1}{p' } \ge  \frac{3}{2q} + \frac{1}{4} > \frac{1}{q}$, i.e. $ p' < q < \sigma $. 

\vspace{0.2cm}
Applying Lemma\,\ref{lemA.1} with $ \Psi = \partial _3\Phi _n, f= u \otimes u \chi _{ Q_0(R)}$ and $ v=u_3$, noting that due to  \eqref{2.4j} and  \eqref{2.4n} $ \Psi $ satisfies  \eqref{A.0a} with $ \alpha =2$, verifying  that $ Q_j(R) \cap Q_0(R)= Q_j(R)$, 
we get from  \eqref{A.1} 
\begin{align}
III_{ 11}  &\le c  \sum_{k=0}^{n} r_k^{ -2 } \|u \otimes u \|_{ L^{ p'}(-r_k^2, 0; L^{ q'}(U_k(R)))}  
 \|u_3\|_{ L^{p}(-r_k^2, 0; L^{ q}(U_k(R)))} 
 \cr
 &\qquad + c\sum_{j=0}^{n}  \sum_{k=j}^{n}    r_k^{ \frac{1}{q'}-2 }  r_{ j}^{ \frac{2}{q'}- \frac{3}{l}} 
  \|u \otimes u\|_{L^{ p'}(-r_k^2, 0;  L^{ l}(U_j(R)))} \|u_3\|_{ L^{p}(-r_k^2, 0; L^{ q}(U_k(R)))}.
   \label{2.50b}
 \end{align}
Recalling the definition of $ Y(r_k)$, noting that $ \frac{3q}{2} > p'$,  applying H\"older's inequality,  together with  \eqref{2.6} we estimate   
for $ k=0, \ldots, n$,
\begin{align*}
r_k ^{ -2}
 & \|u \otimes u \|_{L^{ p'}(-r_k^2, 0;  L^{ q'}(U_k(R)))} \|u_3\|_{ L^{p}(-r_k^2, 0; L^{ q}(U_k(R)))}  
\\
&\le c   r_k^{ -2 } r_k^{ \frac{2q-3p'}{qp'}} \|u \otimes u \|_{L^{ \frac{2q}{3}}(-r_k^2, 0;  L^{ q'}(U_k(R)))} r_k^{ \frac{2}{p}+\frac{3}{q} -1} Y(r_k)
\\
&\le c   r_k^{ -1 }\|u  \|^2_{L^{ \frac{4q}{3}}(-r_k^2, 0;  L^{ 2q'}(U_k(R)))}  Y(r_k)
\\
&\le c   r_k^{ -1 }E_k(R) Y(r_k). 
\end{align*}
Similarly, observing  \eqref{2.50a} and  \eqref{2.6},  we find for $ k,j=0, \ldots, n, k \ge j$, 
\begin{align*}
& r_k^{ \frac{1}{q'}-2 }  r_{ j}^{ \frac{2}{q'}- \frac{3}{l}} 
  \|u \otimes u \|_{L^{ p'}(-r_k^2, 0;  L^{ l}(U_j(R)))} \|u_3\|_{ L^{p}(-r_k^2, 0; L^{ q}(U_k(R)))}  
\\
&\le c   r_k^{ \frac{1}{q'}-2 }  r_{ j}^{ \frac{2}{q'}- \frac{3}{l}} r_k^{ 2\frac{\sigma - p'}{\sigma p'}}
  \|u \|^2_{L^{2 \sigma }(-r_k^2, 0;  L^{ 2l}(U_j(R)))}    r_k^{ \frac{2}{p}- \frac{3}{q} -1} Y(r_k)
  \\
 & = c   r_k^{ \frac{2}{q} - \frac{2}{\sigma } }  r_{ j}^{ \frac{2}{q'}- \frac{3}{l}} 
  \|u \|^2_{L^{2 \sigma }(-r_k^2, 0;  L^{ 2l}(U_j(R)))}     Y(r_k)
  \\
 & \le c   r_k^{ \frac{2}{q} - \frac{2}{\sigma } }  r_{ j}^{ \frac{2}{\sigma } -\frac{2}{q}} r_j^{ -1} E_j(R) Y(r_j).
\end{align*}
Inserting the above estimates into  \eqref{2.50b}, we infer 
\[
III_{ 1} \le c  \sum_{j=0}^{n} r_j^{ -1} E_j(R) Y(r_j) + III_{ 12}  . 
\]
The estimation of $ III_{ 12}$ we postpone after the discussion on $ III_{ 2}$.  By the same argument we have used for 
$ III_1$ we write
\[
III_2= \intl_{-1}^t\intl_{U_0(R)} \pi_0 u\cdot    \Phi _n \nabla (\eta \psi)  dxds + \intl_{-1}^t\intl_{U_0(R)} \pi_h u\cdot    \Phi _n \nabla (\eta \psi)  dxds = III_{ 21}+III_{ 22}. 
\]    
To estimate $ III_{ 21}$ we use Lemma\,\ref{lemA.1} with $ \alpha =1, \Psi =\Phi _n, p=q= \frac{5}{2}$, $ l=\frac{5}{3}, v =u, f= u \otimes u \chi _{ U_0(R)}$ and $ \nabla (\eta \psi) $ in place of $ \eta $. This together with H\"older's inequality and 
 \eqref{2.6} yields
 \begin{align*}
III_{ 21}  &\le c (R-\rho )^{ -1} \sum_{k=0}^{n} r_k^{ -1 } \|u \otimes u \|_
{ L^{ \frac{5}{3}}(-r_k^2, 0; L^{ \frac{5}{3}}(U_k(R)))}  
 \|u\|_{ L^{ \frac{5}{2}}(-r_k^2, 0; L^{ \frac{5}{2}}(U_k(R)))} 
 \cr
 &\qquad + c(R-\rho )^{ -1}\sum_{j=0}^{n}  \sum_{k=j}^{n}    r_k^{- \frac{2}{5}}  r_{ j}^{ - \frac{3}{5}} 
  \|u \otimes u\|_{L^{ \frac{5}{3}}(-r_k^2, 0;  L^{ \frac{5}{3}}(U_j(R)))} 
  \|u\|_{ L^{ \frac{5}{2}}(-r_k^2, 0; L^{ \frac{5}{2}}(U_k(R)))}
\\
&\le c (R-\rho )^{ -1} \sum_{k=0}^{n} r_k^{ - \frac{1}{2} } \|u \|^2_{ L^{ \frac{10}{3}}(Q_k(R))}  
 \|u\|_{ L^{ \frac{20}{3}}(-r_k^2, 0; L^{ \frac{5}{2}}(U_k(R)))} 
 \cr
 &\qquad + c (R-\rho )^{ -1}\sum_{j=0}^{n}  \sum_{k=j}^{n}    r_k^{ \frac{1}{2}- \frac{2}{5}}  r_{ j}^{ - \frac{3}{5}} 
  \|u\|_{L^{ \frac{10}{3}}(Q_j(R))} ^2\|u\|_{ L^{ \frac{20}{3}}(-r_k^2, 0; L^{ \frac{5}{2}}(U_k(R)))}
\\
&\le c (R-\rho )^{ -1} {\cal E}^{ \frac{1}{2}}  \sum_{k=0}^{n} r_k^{ - \frac{1}{2} } E_k(R). 
 \end{align*}
Accordingly
\[
III_2 \le c(R-\rho )^{ -1}  {\cal E}^{ \frac{1}{2}}\sum_{i=0}^{n}r_j^{ - \frac{1}{2}} 
E_j(R)  + III_{ 22}.
\]
It only remains to estimate the sum $III_{ 12}+ III_{ 22}$. In fact, applying integration by parts, and recalling that $ \nabla \cdot u=0$, we calculate 
\begin{align*}
  III_{ 12}+ III_{ 22} &=- \intl_{Q_0(R)} \nabla \pi _h \cdot  u \Phi_n \eta \psi  dxds
\\
&=-  \sum_{k=1}^{n-1} \intl_{A_k(R)} \nabla \pi _h \cdot  u \Phi_n \eta\psi  dxds -  
\intl_{Q_n(R)} \nabla \pi _h \cdot  u \Phi_n \eta \psi  dxds
\\
 &\qquad \qquad - \intl_{A_0(R)} \nabla \pi _h \cdot  u \Phi_n \eta\psi  dxds
\\
&=J_1+J_2+ J_3. 
\end{align*}
Noting that $ 1- \chi_2= 1$ on $ A_0(R)$, using integration by parts, together with  \eqref{2.4j} and  \eqref{hp} we find 
\begin{align*}
  J_3 &=   \intl_{Q_2(R)}  \pi _h \cdot  u \Phi_n \nabla (\eta\psi (1-\chi _2))  dxds 
+  \intl_{Q_2(R)}  \pi _h \cdot  u \partial _3\Phi_n \eta\psi (1-\chi _2)  dxds
\\
& \le c (R- \rho )^{ -1}\|\pi _h\|_{ L^{ \frac{3}{2}}(Q_0(R))} \|u\|_{ L^3(Q_0(R))} 
\\
&\le c (R- \rho )^{ -1} {\cal E}^{ \frac{1}{2}} (\|\pi \|_{ L^{ \frac{3}{2}}(Q_0(R))}   + 
{\cal E}).
\end{align*}
For $ J_1+J_2$ using  \eqref{2.4j} and  \eqref{2.4n}, and applying H\"older's inequality along with  \eqref{2.6}, we get 
\begin{align*}
 J_1+ J_2 &\le c 
\sum_{k=1}^{n} r_k^{ -1}\intl_{Q_k(R)} |\nabla \pi _h|   |u|  dxds 
 \\
&\le c\sum_{k=1}^{n} r_k^{ -1} \|\nabla \pi _h\|_{ L^{ \frac{3}{2}}(-r_k^2, 0; L^{ \infty}(U_k(R)( \frac{R+\rho }{2})))} 
\|u\|_{ L^3(-r_k,0; L^1(U_k(R)))}
\\
&\le c {\cal E}^{ \frac{1}{2}}\sum_{k=1}^{n} r_k^{ \frac{1}{6}}\|\nabla \pi _h\|_{ L^{ \frac{3}{2}}(-1, 0; L^{ \infty}(U_k(R)( \frac{R+\rho }{2})))}.
\end{align*}
Let $ (x, t) \in U_k(R) ( \frac{R+\rho }{2})\times  (-1, 0)$. Since $ \pi _h(\cdot , t)$ is harmonic in $ U_0(R)$ and 
$ \dist (x, \partial U_0(R)) \ge c (R-\rho )$ by using Caccioppli inequality and the mean value property of harmonic functions, we obtain 
\[
|\nabla \pi _h(x, t)| \le c (R-\rho )^{ -1} \frac{1}{|B(R)|}\intl_{B(R)} |\pi _h(\cdot , t)| dx \le c (R-\rho )^{ -3} 
\|\pi _h(\cdot , t)\|_{ L^{ \frac{3 }{2}}(U_0(R))}. 
\] 
This together with  \eqref{hp} shows that 
\begin{align*}
\|\nabla \pi _h\|_{ L^{ \frac{3}{2}}(-1, 0; L^{ \infty}(U_k(R)( \frac{R+\rho }{2})))}
&\le c (R-\rho )^{ -3} \|\pi _h\|_{ L^{ \frac{3}{2}}(Q_0(R))}
 \\
 &\le c (R-\rho )^{ -3} 
(\|\pi \|_{ L^{ \frac{3}{2}}(Q_0(R))} + {\cal E}).    
\end{align*}
Using  the above    estimate, we find 
\[
J_1+J_2 \le c (R-\rho )^{ -3} {\cal E}^{ \frac{1}{2}}
(\|\pi \|_{ L^{ \frac{3}{2}}(Q_0(R))} + {\cal E}).    
\]

Accordingly,
\[
III_{ 12} + III_{ 22} \le c (R-\rho )^{ -3} {\cal E}^{ \frac{1}{2}}
(\|\pi \|_{ L^{ \frac{3}{2}}(Q_0(R))} + {\cal E}).    
\]
Gathering the above estimates,  we deduce 
  \begin{align*}
III & \le c\sum_{j=0}^{n} 
 r_j^{ -1} Y(r_j)E_j(R)+ c (R-\rho )^{ -1} {\cal E}^{ \frac{1 }{2}}\sum_{j=0}^{n}  r_j^{ - \frac{1}{2}} E_j(R)
  \\
 & \qquad + c (R-\rho )^{ -3} {\cal E}^{ \frac{1}{2}}
(\|\pi \|_{ L^{ \frac{3}{2}}(Q_0(R))} + {\cal E}).
    \end{align*}
 On the other hand, from the definition of $ \Phi _n$ for $ n \in \N$ we obtain 
 \begin{align*}
 &\frac{1}{2} \sup_{ t\in (-1,0)}\intl_{ U_0(R)} |u(\cdot , t)|^2 \Phi_n(\cdot, t)  \eta(\cdot , t) \psi dx + 
 \intl_{-1}^0\intl_{ U_0(R)} |\nabla u|^2 \Phi_n  \eta\psi dx ds  
 \\
  & \ge c r_n^{ -1} E_n(\rho ). 
 \end{align*} 
Estimating the left-hand side  of  \eqref{en-bn} from below  using  the above estimate,  and inserting the estimates of  $ I, II$ and $ III$ into the right-hand side of  \eqref{en-bn}, we arrive at 
\begin{align}
  r_n^{ -1} E_n(\rho ) &\le 
c\sum_{j=0}^{n} 
 r_j^{ -1} E_j(R)Y(r_j)+ c (R-\rho )^{ -1} {\cal E}^{ \frac{1 }{2}}\sum_{j=0}^{n}  r_j^{ - \frac{1}{2}} E_j(R)
  \cr
 & + c (R-\rho )^{ -3} {\cal E}^{ \frac{1}{2}}
(\|\pi \|_{ L^{ \frac{3}{2}}(Q_0(R))} + {\cal E}) + c (R-\rho )^{ -2} {\cal E}.   \label{2.22}
\end{align} 
For the second term of the right hand side of \eqref{2.22} we use H\"{o}lder's and Young's inequalities to estimate
\begin{align*}
c (R-\rho )^{ -1} {\cal E}^{ \frac{1 }{2}}\sum_{j=0}^{n}  r_j^{ - \frac{1}{2}} E_j(R)  &\le 
c (R-\rho )^{ -1} {\cal E}^{ \frac{3 }{4}}\sum_{j=0}^{n}  r_j^{ - \frac{5}{8}} E_j(R)^{ \frac{3}{4}} r_j^{ \frac{1}{8}} 
\\
&\le c   (R-\rho )^{ -1} {\cal E}^{ \frac{3 }{4}}\Big(\sum_{j=0}^{n}  r_j^{ - \frac{5}{6}} E_j(R)\Big)^{ \frac{3}{4}}  
\Big(\sum_{j=0}^{n} r_j^{ \frac{1}{2}} \Big)^{ \frac{1}{4}}
\\
&\le \sum_{j=0}^{n}  r_j^{ - \frac{5}{6}} E_j(R) +c   (R-\rho )^{ -4} {\cal E}^{ 3},
\end{align*}
and we deduce from  \eqref{2.22}
\begin{align}
  r_n^{ -1} E_n(\rho ) &\le 
c_0\sum_{j=0}^{n}  r_j^{ -1} E_j(R)Y(r_j)+\sum_{j=0}^{n}  r_j^{ - \frac{5}{6}} E_j(R) 
    + c (R-\rho )^{ -4} C_0,
\label{2.22b}
\end{align} 
where 
\[
C_0= \Big\{1+ {\cal E}^{ 3} + \|\pi \|_{ L^{ \frac{3}{2}}(Q_0(R))}^{ \frac{3}{2}}\Big\}. 
\]

First, by our assumption  \eqref{cond1} it holds  $ Y(R) \rightarrow 0$ as $ R \searrow 0$. Let $ 0< \lambda  <1$ arbitrarily fixed. There exists $ n_0\in \N$ 
such that $ c_0Y(r_{n_0}) + r_{ n_0}^{ \frac{1}{6}} \le \frac{1- 2^{\lambda -1 } }{2}$. Thus,  we get from  \eqref{2.22b}
\begin{equation}
  r_n^{ -1} E_n(\rho )\le  \frac{1- 2^{ \lambda -1 } }{2} \sum_{j=n_0}^{n} r_j^{ -1}E_j(R)
 + c2^{ n_0} ({\cal E}Y(r_0) +{\cal E}) + c (R-\rho )^{ -4} C_0
 \label{2.24a}
 \end{equation} 
 for all $n\ge n_0+1$, while for  $n\le n_0$ we have the following   simple estimate  
\begin{equation}
\label{2.24b}
r_n ^{-1} E_n (\rho) \le  2^{n_0} E_{0} (1)= 2^{n_0} {\cal E}.
 \end{equation} 

Combining \eqref{2.24a} and \eqref{2.24b}, we find
\begin{equation}
  r_n^{ -1} E_n(\rho )\le  \frac{1- 2^{ \lambda -1 } }{2} \sum_{j=1}^{n} r_j^{ -1}E_j(R)
 + c2^{ n_0} ({\cal E}Y(r_0) +{\cal E}) + c (R-\rho )^{ -4} C_0
 \label{2.24}
 \end{equation} 
for all $n\in \Bbb N$.  Given $ N\in \N$, multiplying  \eqref{2.24} by $ 2^{- n(1-\lambda ) }= r_n^{ 1-\lambda }, $ and  summing it  from $n=1$ to $ n=N$, we obtain 
\begin{align*}
  &\sum_{ n=1}^{N}  r^{ -\lambda    }_n E_n(\rho ) 
  \\
  &\le  
  \frac{1- 2^{ \lambda -1 } }{2} \sum_{ n=1}^{N}\sum_{j=1}^{n} r_n ^{ 1-\lambda  }r_j^{ -1} E_j(R)+ 
  c\Big\{2^{ n_0} ({\cal E}Y(r_0) +{\cal E}) + (R-\rho )^{ -4} C_0\Big\} \sum_{ n=1}^{N} r_n^{ 1-\lambda  }
 \\
 &\le   \frac{1- 2^{ \lambda -1 } }{2} \sum_{j=1}^{N}\sum_{n=j}^{N} r_n ^{1-\lambda  } r_j^{ -1} E_j(R) +
  c   \Big\{2^{ n_0} ({\cal E}Y(r_0) +{\cal E}) +  (R-\rho )^{ -4} C_0\Big\}  
\\
 &\le  \frac{1 }{2}   \sum_{j=1}^{N}r_j^{ -\lambda  } E_j(R) +   c  \Big\{2^{ n_0} ({\cal E}Y(r_0) +{\cal E}) +  (R-\rho )^{ -4} C_0\Big\}. 
\end{align*}
Applying the algebraic lemma \cite[V. Lemma\,3.1]{gia}, we arrive at  
\begin{equation}
   \sum_{n=1}^{N}  r^{ -\lambda    }_n E_n\Big(\frac12\Big)  \le c \Big\{2^{ n_0} ({\cal E}Y(r_0) +
   {\cal E}) +  C_0\Big\}. 
\label{2.27}
 \end{equation} 
In particular, we get for all $\lambda, r \in (0,1)$ the Morrey-type estimate 
\begin{equation}
r^{ -\lambda  } \|u\|^2_{ L^\infty(-r^2, 0; L^2(B'(r)\times  (-r, r)))} + 
r^{ -\lambda  } \|\nabla u\|^2_{ L^2(B'(r )\times  (-r,r)\times  (-r^2, 0)) )} \le C, 
\label{2.28}
 \end{equation} 
 where $ C$ depends on $ \lambda $. This completes proof of the first statement of the theorem.  
 
 \vspace{0.3cm}
{\it Proof of the second statement}.  Now, we assume the second condition  \eqref{cond2} is fulfilled for some $ \alpha \in (1, +\infty)$.  Let $ \varepsilon _0>0$ be sufficiently small to be specified below. Observing  \eqref{cond2}, 
there exists $ n_0\in \N$ such that $ j^{ \alpha } (c_0Y(r_{j}) + r_j^{ \frac{1}{6}}) \le \varepsilon _0 $ for all $ j \ge n_0$.  
Thus,   \eqref{2.22} together with H\"older's inequality and Young's inequality yields
 \begin{align}
r_n^{ -1} E_n(\rho ) &\le   \varepsilon_0 \sum_{j=1}^{n} j^{ -\alpha } r_j^{ -1} E_j(R) +(R-\rho )^{ -4}C_1,
\label{2.29}
\end{align}
where $ C_1=c  \Big\{2^{ n_0} ({\cal E}Y(r_0) +{\cal E}) +   C_0\Big\}$. 
Given $ N\in \N$, we multiply $  \eqref{2.29}$ by $ n^{ -\alpha }$ and summing both sides of the resultant 
inequality from $ n=1$ to $ N$. This leads to 
\begin{align*}
 \sum_{n=0}^{N} n^{ -\alpha }r_n^{ -1} E_n(\rho ) &\le   \varepsilon_0 
 \sum_{n=0}^{N} \sum_{j=0}^{n}n^{ -\alpha } j^{ -\alpha } r_j^{ -1} E_j(R) +c(R-\rho )^{ -4}C_1  \sum_{n=0}^{N} n^{ -\alpha }
\\
&=\varepsilon _0\sum_{j=0}^{N} \sum_{n=j}^{n} n^{ -\alpha }j^{ -\alpha } r_j^{ -1} E_j(R) +c(R-\rho )^{ -4}C_1
\\
&  \le c_{ \alpha }\varepsilon _0\sum_{j=0}^{N} j^{ -2\alpha +1}r_j^{ -1} E_j(R) +c(R-\rho )^{ -4}C_1
\\
&\le c_{ \alpha }\varepsilon _0\sum_{j=0}^{N} j^{ -\alpha }r_j^{ -1} E_j(R) +c(R-\rho )^{ -4}C_1.
\end{align*}
Taking $ \varepsilon _0 = \frac{1}{2c_\alpha }$ and once more using theiteration lemma    \cite[V. Lemma\,3.1]{gia}, we deduce that 
\begin{equation}
  \sum_{j=0}^{N} j^{ -\alpha }r_j^{ -1} E_j\Big( \frac{3}{4}\Big) \le cC_1.
\label{2.29b}
 \end{equation}   
Thus,  \eqref{2.29} with $ \rho = \frac{1}{2}$ and $ R= \frac{3}{4}$ gives 
\begin{equation}
r_n^{ -1} E_n \Big(\frac{1}{2}\Big) \le cC_1\quad \forall n\in \N. 
\label{2.30}
 \end{equation} 
 As a consequence of  \eqref{2.30} we get 
\begin{equation}
K_0:=\sup_{0< r \le  1}r^{ -1} E(r) < +\infty,   
\label{2.31c}
 \end{equation} 
where 
\[
E(r)=E^{ (1)}(r)+E^{ (2)}(r):=  \esssup_{t\in (-r^2, 0)}\intl_{B(r)} |u(t)|^2 dx + \intl_{Q(r)} |\nabla u|^2 dx ds. 
\]
Using H\"older's inequality along with  \eqref{2.6},  we find 
\begin{equation}
r^{ -2}\|u\|_{ L^3(Q(r))}^3 \le  c r^{ -\frac{3}{2}} \|u\|_{ L^4(-r^2, 0; L^3(B(r)))}^3 \le c r^{- \frac{3}{2}} E(r)^{ \frac{3}{2}}
\le c  K_0^{ \frac{3}{2}}.  
\label{2.31a}
 \end{equation}

{\bf Estimation of the pressure:}  Similarly to   \cite[Lemma\,2.9]{cw6} we get the following 
\begin{lem}
Let $s\in (1, +\infty)$, and $ f\in L^{ s}(Q(1))$   satisfy
\begin{equation}
\intl_{Q(r)} |f| ^s dxds \le K_0^s r^{ \lambda }\quad  \forall  r_0 \le  r \le  1
\label{2.44}
 \end{equation} 
  for some $ 0< \lambda  < 5$ and $ 0< r_0 <1$. 
Let $ \pi \in L^s(Q(1))$ solve  $\Delta  \pi = \partial _i \partial _j f_{ ij} $ in the sense of distributions. 
Then 
\begin{equation}
\intl_{Q(r)} |\pi - \pi _{ B(r)}| ^s dxds \le c r^{ \lambda } \intl_{Q(r)} |\pi - \pi _{ B(r)}| ^s dxds +cK_0^s r^{ \lambda }\quad  \forall  r_0 \le  r \le  1,
\label{2.45}
 \end{equation}
 where $ c=\const>0$ depends only on $ \lambda$ and $ s$. 
\end{lem}

{\bf Proof}: Let $ r_0 \le r \le 1$ and $ \theta \in (0,1)$. We write 
\begin{equation}\label{2.45a}
 \pi(t) - \pi(t) _{ B(r)} = \pi _0(t) + \pi _h(t),
  \end{equation}
   where 
$ \pi _0 \in A^s(B(r)) = \{v= \Delta q\,| \,q\in W^{2,\,s}_0(B(r)) \}$, and $\Delta \pi_h =0$ on $B(r)$. Thus, in view of \cite[Lemma\,2.8]{cw6} we estimate 
\[
\|\pi_0 (t)\|^s_{ L^s(B(r))} \le c\|f(t)\|^s_{ L^s(B(r))}.
\]
Furthermore, using the mean value property together with the Caccioppoli inequality for harmonic functions, we obtain 
\begin{align}\label{2.45b}
\|\pi _h(t) - (\pi _h(t))_{ B(\theta r)}\|_{ L^s(B(\theta r))}^s &\le c (\theta r)^5 \|\nabla \pi _h(t)\|_{ L^\infty(B(\theta r))}^s \cr
&\le c \theta ^5 \| \pi _h(t)\|_{L^s(B( r))}^s \le   c \theta ^5 \| \pi(t)- \pi (t)_{ B(r)}\|_{L^s(B( r))}^s.  
\end{align} 
By the  triangle inequality  together with \eqref{2.45a} and \eqref{2.45b} we get
\begin{align*}
  &\| \pi(t)- \pi(t) _{ B(\theta r)}\|_{L^s(B(\theta  r))}^s  
\\
&\le\|\pi _h(t) - (\pi _h(t))_{ B(\theta r)}\|_{ L^s(B(\theta r))}^s+ \|\pi _0(t) - (\pi _0(t))_{ B(\theta r)}\|_{ L^s(B(\theta r))}^s
\\
&\le   c \theta ^5 \| \pi(t)- \pi (t)_{ B(r)}\|_{L^s(B( r))}^s+c\|f(t)\|^s_{ L^s(B(r))}.
\end{align*}
Integrating both sides over $ (-(\theta r)^2, 0)$ with respect to $ t$ and applying  \eqref{2.44}, we see that 
\begin{align*}
  &\| \pi- \pi _{ B(\theta r)}\|_{L^s(Q(\theta  r))}^s  \le   c \theta ^5 \| \pi- \pi _{ B(r)}\|_{L^s(Q( r))}^s+c K_0^s r^{ \lambda }.  
\end{align*}
 Given $ \lambda < \mu  < 5$, we may choose $ \theta $ such that $ c \theta ^{ 5-\mu } \le 1$. This together with a standard iteration yields 
\[
\| \pi- \pi _{ B(r)}\|_{L^s(Q(r))}^s  \le c r^{ \lambda } \| \pi(t)\|_{L^s(Q(1))}^s + c r^{ \lambda } K_0^s\quad  \forall 
r_0 \le r \le 1. 
\]
Whence the claim.  \hfill \Beweisende

\vspace{0.3cm}
Applying the above lemma for $ s = \frac{3}{2}, f= u \otimes u, \lambda =2$ and taking into account \eqref{2.31a}, we get 
\begin{equation}
  r^{ -2}\|\pi - \pi _{ B(r)}\|_{ L^{ \frac{3}{2}}(Q(r))}^{ \frac{3}{2}} \le c  K_0^{ \frac{3}{2}} \quad  \forall 0< r< 1. 
\label{2.31b}
 \end{equation} 
 
{\bf  Proof that  $ (0,0)$ is regular point of $ u$ via indirect argument}.  Before starting  the proof  we 
recall the $ \varepsilon $-regularity condition proved in \cite{wol8}. 
There is an absolute number $ \varepsilon _0>0$ with the following property. If $ u\in V^{1, 2}(Q(1))$ is suitable weak solution of  \eqref{nse} and there exists $ 0<r \le 1$ such that 
\begin{equation}
r^{ -2} \intl_{Q(r)} |u|^3 dxdt \le \varepsilon _0,  
\label{eps}
 \end{equation} 
then $ u\in L^\infty(Q( \frac{r}{2}))$, in particular $ (0,0)$ is a regular point of $ u$. 
(In fact the above criterion is proved  in \cite{wol8} for more general notion of {\em local suitable weak solution}.) 
\vspace{0.2cm}
 Now, assume $ (0,0)$ is not a regular point of $ u$. In view of the $ \varepsilon $-regularity criterion  we stated above 
 \begin{equation}
  r^{ -2} \|u\|_{ L^3(Q(r))}^3 > \varepsilon _0\quad \forall 0< r \le 1, 
 \label{2.46}
  \end{equation} 
  where $ \varepsilon _0 >0$ denotes the absolute number in  \eqref{eps}.  
  Define,  with  $ r_k= 2^{ -k}, k\in \N$,
\begin{align*}
v_k(x,t ) &= r_k u(r_k x,  r_k^2 t ),\quad     
\\
\pi _k(x,t ) &= r^2_k \Big(\pi (r_k x,  r_k^2 t ) -\pi(r_k^2t  )_{ B(r_k)}\Big),\quad    (x, t )\in Q(1). 
\end{align*}  
Then $ (v_k,\pi _k)\in V^{ 1,2}(Q(1))\times  L^{ \frac{3}{2}}(Q(1)), k\in \N, $ is a suitable weak solution to the Navier-Stokes 
equations in $ Q(1)$. In view of  \eqref{2.31c} and  \eqref{2.31b}, using a standard scaling argument, we see that   $ \{(v_k,\pi _k)\}$ is bounded in $ V^{ 1,2}(Q(1))\times  L^{ \frac{3}{2}}(Q(1))$. Furthermore, 
 \eqref{2.46} turns into 
 \begin{align}
r^{ -2}\|v_k\|_{ L^3(Q(r))}^3 & > \varepsilon _0\quad  \forall 0< r< 1,
   \label{2.46a} 
 \end{align}

Eventually passing to a subsequence, we get $ (\overline{v},\overline{\pi} )\in V^{ 1,2}(Q(1))\times  L^{ \frac{3}{2}}(Q(1))$ such that 
\begin{align}
  \nabla v_k &\rightarrow \nabla \overline{v} \quad   \text{weakly in}\quad L^2(Q(1)) \quad  \text{as}\quad k \rightarrow +\infty, 
   \label{2.48}
\\
  v_k &\rightarrow  \overline{v} \quad   \text{weakly$-\ast$ in}\quad L^\infty(-1,0; L^2(B(1))) \quad  \text{as}\quad k \rightarrow +\infty, 
\label{2.49}
\\
\pi _k &\rightarrow \overline{\pi } \quad   \text{weakly in}\quad L^{ \frac{3}{2}}(Q(1)) \quad  \text{as}\quad k \rightarrow +\infty. 
\label{2.50}
\end{align}
Furthermore, using Lions-Aubin's Lemma, we get  
\begin{align}
   v_k &\rightarrow \overline{v} \quad   \text{strongly in}\quad L^3(Q(1)) \quad  \text{as}\quad k \rightarrow +\infty. 
\label{2.51}
\end{align}
In particular, $ (\overline{v}, \overline{\pi })$ is a suitable weak solution to the Navier-Stokes equations. By the aid of  
\eqref{2.51} in \eqref{2.46a}  letting $ k \rightarrow +\infty$, we obtain
\begin{align}
r^{ -2}\|\overline{v}\|_{ L^3(Q(r ))}^3  & \ge \varepsilon _0\quad  \forall 0< r< 1.  
 \label{2.48b}
 \end{align} 
On the other hand, from  \eqref{cond2} (or even from  \eqref{cond1}) we deduce that 
\[
\|v_{k, 3}\|_{ L^{ p}(-1, 0; L^{ q}(B(1)))} = r_k^{ 1- \frac{2}{p}- \frac{3}{q}}\|u_{3}\|_{ L^{ p}(-r_k^2, 0; L^{ q}(B(r_k)))}
\rightarrow 0 \quad  \text{as}\quad  k \rightarrow +\infty.  
\]
Accordingly, $ \overline{v}_3 \equiv 0$. By the localized version of  the one velocity component criterion in  \cite{nenope} 
(see also \cite[Section\,1.4]{nest}) we find  that $\overline{v}\in L^\infty (Q( \frac{1}{2}))$,  and the left hand side of  \eqref{2.48b}  goes to zero as $r\to 0$, which is a contradiction.  Consequently the assumption is not true and therefore $ (0,0)$ is regular point. 
 \hfill \Beweisende

 \hspace{0.5cm}
$$\mbox{\bf Acknowledgements}$$
Chae was partially supported by NRF grant 2016R1A2B3011647, while Wolf has been supported 
supported by NRF grant 2017R1E1A1A01074536.
The authors declare that they have no conflict of interest.

\appendix
%%% ----------------------------------------------------------------------
%       SECTION 
%%% ----------------------------------------------------------------------
\section{Appendix}
\label{sec:-A}
\setcounter{secnum}{\value{section} \setcounter{equation}{0}
\renewcommand{\theequation}{\mbox{A.\arabic{equation}}}}

The aim of this appendix is to provide an estimate which will be used various times  for the estimation of integrals 
involving the pressure during the proof of our main result. 

We start our discussion to define the following singular integral operator, which is nrcessary  for the decomposition of the pressure. 
Throughout this appendix, let $ 0< R \le 1$ be fixed.  
Let $ Q_0(R) = U_0(R) \times  (-1,0)$, where $ U_0(R) = B'(R)\times  (-1,1)$. Given $ f_{ ij}\in L^p(Q_0(R)), 1<p< +\infty, i,j=1,2,3,$ we define 
\[
\mathscr{J}(f)(x,t) = P.V. \intl_{ \R^{3} } K(x-y) :f(y,t) \chi _{ U_0(R)}(y)dy,\quad  (x,t)\in \R^{3}\times  (-1,0), 
\]
with the Calder\'on-Zygmund kernel  $ K_{ ij}= \partial _{ i}\partial _j N, i,j=1, 2,3,$  where 
\[
N(x) = \frac{1}{4\pi |x|},\quad  x\in \R^{3} \setminus \{0\}.  
\]
Clearly by virtue  of Calder\'on-Zygmund inequality, $ \mathscr{J}: L^p(Q_0(R)) \rightarrow L^p(\R^{3} \times  (-1,0))$ defines a bounded linear operator. In particular, 
\begin{equation}
 \|\mathscr{J}(f)\|_{ L^p(\R^{3} )} \le c \|f\|_{ L^p(U_0(R))}. 
\label{A.6}
 \end{equation} 
Furthermore, 
setting $ \pi _0= \mathscr{J}(f)$,  it holds 
\begin{equation}
 -\Delta  \pi _0 = \nabla \cdot \nabla \cdot f\quad   \text{in}\quad  Q_0(R)
\label{A.0}
 \end{equation} 
in the sense of distributions. 
As in Section\,2 we use the following  notation 
\[
r_j= 2^{ -j},\quad  U_j(R) = B'(1)\times  (-r_j, r_j),\quad Q_j(R) = U_j(R) \times  (-r_j^2, 0)\quad  j\in \N_0. 
\]

\begin{lem}
 \label{lemA.1}
 Let $ n\in \N$. Let $ \Psi \in C^\infty(\R^3\times  (-\infty, 0))$ such that for constants $ \alpha>0, c>0 $ and $ C>0$ 
 it holds 
  \begin{equation}
 \begin{cases}
 c r_j^{ -\alpha } \le \Psi (x, t) \le Cr_j^{ -\alpha } \quad  \forall (x, t)\in A_j =Q_j(R) \setminus Q_{ j+1}(R),\quad \forall j=0, \ldots, n-1
 \\[0.3cm]
 c r_n^{ -\alpha } \le \Psi (x, t) \le Cr_j^{ -\alpha } \quad  \forall (x, t)\in Q_n(R). 
 \end{cases} 
 \label{A.0a}
  \end{equation} 
 Let $ 1 < p,q <+\infty, 1 < l \le  q'$. Let $ v\in L^{ p}(-1, 0; L^q(U_0(R)))$, and 
$ f\in L^{ p'}(-1, 0; L^{ q'}(U_0(R)))$,   $ 1< m, l < +\infty$. 
Then, setting $ \pi _0 = \mathscr{J}(f)$ it holds 
\begin{align}
 & \intl_{-1}^{ t}\intl_{U_0(R)} \pi _0 v \Psi  \eta dxds  
\cr
 &\le c \sup|\eta | \sum_{k=0}^{n} r_k^{ -\alpha } \|f\|_{ L^{ p'}(-r_k^2, 0; L^{ q'}(U_k(R)))}  
 \|v\|_{ L^{p}(-r_k^2, 0; L^{ q}(U_k(R)))} 
 \cr
 &\qquad + c\sup|\eta |\sum_{j=0}^{n}  \sum_{k=j}^{n}    r_k^{ \frac{1}{q'}-\alpha }  r_{ j}^{ \frac{2}{q'}- \frac{3}{l}} 
  \|f\|_{L^{ p'}(-r_k^2, 0;  L^{ l}(U_j(R)))}
 \|v\|_{ L^{p}(-r_k^2, 0; L^{ q}(U_k(R)))}, 
\label{A.1}
\end{align}
where $ \eta \in  C^\infty_c(U_0(R)\times  (-1, 0])$ stands for a cut off function.  The constant in  \eqref{A.1} depends only on 
$ p, q$ and  $l$.

\end{lem}

{\bf Proof:} Let $ f\in L^p(Q_0(R))$. Set $ \pi _0 = \mathscr{J}(f)$. For $ j\in \N_0$ let $ \chi _j\in 
C^\infty_c(U_j(R) \times  (-r_j^2, 0] )$ with $ \chi _j=1$  on $ Q_{ j+1}(R)$ such that $ 0 \le \chi _j \le 1$, 
and $ |\partial _3\chi_j| \le c r_j^{ -1}$, and $ |\nabla ' \chi_j| \le c R^{ -1}$.	
We set
\[
\phi _j=\begin{cases}
1- \chi _0\quad   &\text{if}\quad  j=0,
\\[0.3cm]
\chi _{ j}- \chi _{ j+1} &\text{if}\quad  j=1, \ldots, n-1,
\\[0.3cm]
\chi _n\quad   &\text{if}\quad  j=n. 
\end{cases}
\]
We have $  \dsum_{j=0}^n\phi _j= 1- \chi _0 + \chi_0- \chi _1 + \ldots+\chi_{ n-1}- \chi _n + \chi _n=1$.  
Accordingly,  $ f =  \dsum_{j=0}^{n} f \phi _j$,
and therefore it holds 
\[
 \pi _0 = \mathscr{J}(f) =  \sum_{j=0}^{n} \mathscr{J}(\phi _j f) = 
\sum_{j=0}^{n} \pi _{ 0,j}. 
\]

This yields
\begin{align*}
&\intl_{-0}^{ t}\intl_{U_0(R)} \pi _0 v\Psi\eta dxds 
=  \sum_{k=0}^{n}\intl_{-1}^{ t}\intl_{U_0(R)} \pi _0 v \Psi \phi _k\eta dxds 
\\
&\quad =   \sum_{j=0}^{n}\sum_{k=0}^{n} \intl_{-1}^{ t}\intl_{U_0(R)} \pi_{ 0,j} v \Psi  \phi _k \eta dxds 
=   \sum_{k=0}^{n}\sum_{j=k}^{n}\intl_{-1}^{ t}\intl_{U_0(R)} \pi_{ 0,j}   v\Psi\phi _k \eta dxds   
\\
&\quad \qquad + \sum_{j=0}^{n}  \sum_{k=j+1}^{n}\intl_{-1}^{ t}\intl_{U_0(R)} \pi_{ 0,j} v \Psi  \phi _k \eta dxds  
= I+II. 
\end{align*}
First, we calculate  
\begin{align*}
  I = \sum_{k=0}^{n}\intl_{-1}^{ t}\intl_{U_0(R)} \Pi_{ 0,k} v\Psi \phi _k\eta dxds,  
\end{align*}
where 
\[
\Pi _{ 0,k}= \begin{cases}
\pi _0\quad  &  \text{if}\quad  k=0,
\\[0.3cm]
  \mathscr{J} (\phi_k f)\quad   &\text{if}\quad  k=1, \ldots, n. 
\end{cases}
\]
Observing  \eqref{A.0a}, applying H\"older's inequality  along with  \eqref{A.6}, we get 
\begin{align*}
I &\le c \sup|\eta | \sum_{k=0}^{n} r_k^{ -\alpha } \|\Pi _{ 0,k}\|_{ L^{ p'}(-r_k^2, 0; L^{ q'}(U_k(R)))}  
 \|v\|_{ L^{ p}(-r_k^2, 0; L^{ q}(U_k(R)))}  
\\
&\le c  \sup|\eta |\sum_{k=0}^{n} r_k^{ -\alpha } \|f\|_{ L^{ p'}(-r_k^2, 0; L^{ q'}(U_k(R)))}  
 \|v\|_{ L^{p}(-r_k^2, 0; L^{q}(U_k(R)))}.   
\end{align*}
For the second integral we find 
\begin{align*}
II &=   \sum_{j=n-2}^{n} \sum_{k=j}^{n} \intl_{-1}^{ t}\intl_{U_0(R)} \Pi_{ 0,j} v\Psi \phi_{ k} \eta dxds  
\\
&\qquad +\sum_{j=0}^{n-3} \sum_{k=j}^{j+3} \intl_{-1}^{ t}\intl_{U_0(R)} \Pi_{ 0,j} v\Psi \phi _{k} \eta  dxds  
\\
&\qquad +  \sum_{j=0}^{n-3}  \sum_{k=j+4}^{n}\intl_{-1}^{ t}\intl_{U_0(R)} \Pi_{ 0,j} v\Psi  \phi _k \eta dxds  = II_1+ II_2+ II_3. 
\end{align*}
Arguing as above, observing  \eqref{A.0a} and applying H\"older's inequality and  \eqref{A.6}, we see that 
\[
II_1+II_2 \le c \sup|\eta | \sum_{k=1}^{n} r_k^{ -\alpha } \|f\|_{ L^{ p'}(-r_k^2, 0; L^{ q'}(U_k(R)))}  
 \|v\|_{ L^{ p}(-r_k^2, 0; L^{ q}(U_k(R)))}.   
\]
It remains to estimate $ II_3$. We calculate 
\begin{equation}
 II_3 =  \sum_{j=0}^{n-2}  \sum_{k=j+4}^{n}\intl_{-1}^{ t}\intl_{U_0(R)} \Pi_{ 0,j} v \Psi \phi _k   dxds =  \sum_{j=0}^{n-2}  \sum_{k=j+4}^{n} J_{ jk} .   
\label{A.7}
 \end{equation} 
Let   $ j+4 \le k \le n$ be fixed. 
Applying H\"older's inequality together with  \eqref{A.0a}, we find 
\begin{align*}
J_{ jk} &\le  c\sup|\eta | r_k^{ -\alpha } \|\Pi _{ 0, j}\|_{ L^{ p'}(-r_k^2, 0; L^{ q'}(U_k(R)))}  
 \|v\|_{ L^{p}(-r_k^2, 0; L^{ q}(U_k(R)))}.
\end{align*}

From  the definition of $ \Pi _{ 0, j}$ it follows that $ \Delta \pi _{ 0, j} = \nabla \cdot \nabla \cdot (f\phi _j)$ 
in the sense of distributions. 
Since $ \supp(\phi _j)  \subset Q_{ j}(R) \setminus Q_{ j+2}(R)$ the function $ \pi _{ 0, j}$ is harmonic in $ \R^{2}\times (-r_{ j+2}, r_{ j+2})\times  (-r_{ j+2}^2, 0)$. Applying Lemma\,\ref{lemA.2} below for $ h=\pi _{ 0, j}, r=r_{ j+2} $ and $ \rho = r_k$, 
we get for almost all $ s\in (-r_k^2, 0)$
\begin{align*}
 \|\Pi _{ 0, j}(s)\|_{ L^{ q'}(U_k(R))}=\|\pi _{ 0, j}(s)\|_{ L^{ q'}(B'(R)\times  (-r_k, r_k))}   \le c\sup|\eta | r_k^{ \frac{1}{q'}}  r_{ j+2}^{ \frac{2}{q'}- \frac{3}{l}} 
  \|\Pi _{ 0, j}(s)\|_{ L^{ l}(\R^{3} )}. 
\end{align*}
Taking the $ L^{ p'}$ norm with respect to $ s$, and employing  \eqref{A.6}, we find 
\begin{align*}
 \|\Pi _{ 0, j}\|_{L^{ p'}(-r_k^2,0;  L^{ q'}(U_k(R)))}   \le c\sup|\eta | r_k^{ \frac{1}{q'}}  r_{ j+2}^{ \frac{2}{q'}- \frac{3}{l}} 
  \|f\|_{L^{ p'}(-r_k^2, 0;  L^{ l}(U_j(R)))}. 
\end{align*}
Accordingly, 
 \begin{align*}
J_{ jk} &\le  c \sup|\eta | r_k^{ \frac{1}{q'}-\alpha }  r_{ j+2}^{ \frac{2}{q'}- \frac{3}{l}} 
  \|f\|_{L^{ p'}(-r_k^2, 0;  L^{ l}(U_j(R)))}
 \|v\|_{ L^{p}(-r_k^2, 0; L^{ q}(U_k(R)))}.
\end{align*}
Inserting this inequality into  \eqref{A.7}, we arrive at 
\begin{align*}
II_3 &=  c\sup|\eta |\sum_{j=1}^{n-2}  \sum_{k=j+3}^{n}    r_k^{ \frac{1}{q'}-\alpha }  r_{ j+2}^{ \frac{2}{q'}- \frac{3}{l}} 
  \|f\|_{L^{ p'}(-r_k^2, 0;  L^{ l}(U_j(R)))}
 \|v\|_{ L^{p}(-r_k^2, 0; L^{ q}(U_k(R)))}
 \\
 &\le   c\sup|\eta |\sum_{j=0}^{n}  \sum_{k=j}^{n}    r_k^{ \frac{1}{q'}-\alpha }  r_{ j}^{ \frac{2}{q'}- \frac{3}{l}} 
  \|f\|_{L^{ p'}(-r_k^2, 0;  L^{ l}(U_j(R)))}
 \|v\|_{ L^{p}(-r_k^2, 0; L^{ q}(U_k(R)))}.
\end{align*}  
Combining the above estimates,  we get the claim.  \hfill \Beweisende

\begin{lem}
 \label{lemA.2}
Let $ 0< r \le  R <+\infty$. Let $ h: B'(2R)\times (-r, r) \rightarrow \R$ be harmonic. Then for all $ 0< \rho \le \frac{r}{4} $ and   
$ 1 \le  l \le p \le  +\infty$ we get
\begin{equation}
 \|h\|_{ L^p(B'(R)\times  (-\rho , \rho))}^p \le c\rho  r^{2 -3\frac{p}{l}}
  \|h\|_{ L^l(B'(2R)\times  (-r, r))}^p, 
  \label{A.5}
 \end{equation} 
where $ c$ stands for a positive constant depending only on $ p$ and $ l$. 
\end{lem} 
 
{\bf  Proof}: Let $ k\in \N, k \ge 2$. Set $ \rho _k = 2^{ -k} r $.  Since $ B'(2R)\times  (- r, r)$ is a non isotropic cylinder, in order to apply the mean value 
property of harmonic functions we use a covering argument.   
We may choose a finite  family of points  $ \{x'_{\nu }\}$ in $ B'(R)$ such that  
$ \{B'(x_{\nu }', r/4)\}$ is a covering of $ \overline{B'(R)}$,  and it holds 
\begin{equation}
 \sum_{\nu } \chi _{ B'(x'_{ \nu }, r)} \le N,\quad  |x_{\nu }- x_{ \mu }| \ge \frac{r}{4}\quad  \forall \nu \neq \mu,   
\label{A.4}
 \end{equation} 
where $ N$ stands for an absolute number. Setting $ x_{\nu }= (x'_{ \nu }, 0)$, we see that 
$ B'(x'_{\nu}, r/4 )\times  (- r/4, r/4)  \subset B(x_{\nu }, r/2)$. With this notation we have 
\begin{align}
  \|h\|_{ L^p(B'(R)\times  (-\rho _k, \rho _k))}^p  
&\le \sum_{\nu }\|h\|_{ L^p(B'(x'_{ \nu }, r/4)\times  (- \rho _k, \rho _k))}^p  
\cr
&\le c r^2 \rho _k\sum_{\nu }\|h\|_{ L^\infty(B'(x'_{ \nu }, r/4)\times  (- r/4, r/4))}^p
\cr
&\le c r^2 \rho _k\sum_{\nu }\|h\|_{ L^\infty(B(x_{ \nu }, r/2))}^p.   
 \label{A.3a}
\end{align}
Since $ h$ is harmonic, using the mean value 
property, we find  
\begin{align*}
\|h\|^p_{ L^\infty(B(x_{\nu}, r/2 ))} 
&\le c r^{ -\frac{3p}{l}}   \|h\|^p_{ L^l(B(x_{\nu}, r ))}  
\\
&\le c  r^{ -\frac{3p}{l}}  \|h\|^{ p-l}_{ L^l( B'(2R)\times  (-r,r))} \|h\|^l_{ L^l(B'(x'_{\nu}, r )\times  (-r,r))}
\\
&= c  r^{ -\frac{3p}{l}}  \|h\|^{ p-l}_{ L^l( B'(2R)\times  (-r,r))} \intl_{-r}^r
\intl_{B'(2R)} |h|^l \chi _{ B'(x'_{\nu}, r )} dx'dx_3.  
\end{align*}
Taking the sum over $ \nu $ and using  \eqref{A.4}, we obtain 
\begin{align}
\sum_{\nu }\|h\|_{ L^\infty(B(x_{ \nu }, r/2))}^p
\le   cN  r^{ -\frac{3p}{l}}  \|h\|^{ p}_{ L^l( B'(2R)\times  (-r,r))}. 
\label{A.3}
\end{align}
Combing  \eqref{A.3a} and  \eqref{A.3},  we get 
\begin{equation}
 \|h\|_{ L^p(B'(R)\times  (-\rho_k , \rho_k))}^p \le c\rho_k  r^{2 -3\frac{p}{l}}
  \|h\|_{ L^l(B'(2R)\times  (-r, r))}^p.  
\label{A.7a}
 \end{equation} 
 Let $ 0< \rho \le  \frac{r}{4}$. Then there exists a unique integer  $ k\ge 2$ such that $ \rho_{ k+1} < \rho \le \rho_k$, 
 Thus,  \eqref{A.7a} implies 
 \begin{align*}
 \hspace{-0.5cm}\|h\|_{ L^p(B'(R)\times  (-\rho , \rho))}^p&\le c\rho_k  r^{2 -3\frac{p}{l}}
  \|h\|_{ L^l(B'(2R)\times  (-r, r))}^p 
 \le  2c\rho  r^{2 -3\frac{p}{l}}
  \|h\|_{ L^l(B'(2R)\times  (-r, r))}^p.    
 \end{align*}
 Whence,  \eqref{A.5}.  \hfill \Beweisende

 \end{document}